\documentclass[12pt,letterpaper]{amsart}
\usepackage[top=2in, bottom=1.5in, left=1in, right=1in]{geometry}
\usepackage{amsmath, amsthm, amssymb, graphicx, url, pstricks, color, colortbl, stmaryrd}

\long\def\forget#1\forgotten{}
\newcommand{\nc}{\newcommand}
\nc{\seq}[1]{#1_1,\allowbreak#1_2,\allowbreak\dots}
\nc{\Set}{\op{Set}}
\nc{\rd}[1]{\color{red}{#1}} \nc{\gr}[1]{\color{green}{#1}}
\nc{\board}[9]{\begin{tabular}{|c|c|c|c|c|c|c|c|c|} \hline #1 & #2 & #3 & #4 & #5 & #6 & #7 & #8 & #9\\ \hline \end{tabular}}

\nc{\sone}{\mathsf{S}_{1}}\nc{\sfin}{\mathsf{S}_\mathrm{fin}}\nc{\ufin}{\mathsf{U}_\mathrm{fin}}
\nc{\Split}{\mathsf{Split}}\nc{\gone}{\mathsf{G}_{1}}\nc{\goab}{\gone(\mathcal{A},\scrB)}
\nc{\gfin}{\mathsf{G}_\mathrm{fin}}
\nc{\fd}{\mathfrak{d}}
\nc{\Impl}{\Rightarrow}
\nc{\Op}{\mathrm{O}}\nc{\Om}{\Omega}\nc{\Ga}{\Gamma}\nc{\Tau}{\mathrm{T}}
\nc{\Asc}{\mathrm{Asc}}
\nc{\ba}{\mathbf{a}}\nc{\bb}{\mathbf{b}}\nc{\bc}{\mathbf{c}}
\nc{\sonedown}{\mathsf{D}}
\nc{\as}{\sub^*}
\nc{\my}[1]{{\sf (#1)}\marginpar{**}}
\nc{\la}{\langle}\nc{\ra}{\rangle}
\nc{\Un}{\bigcup}
\nc{\Inf}{[\bbN]^\infty}
\nc{\Mono}{\op{Mono}}
\nc{\smx}[1]{\bigl(\begin{smallmatrix}#1\end{smallmatrix}\bigr)}
\nc\bbZ{\mathbb{Z}}\nc\bbQ{\mathbb{Q}}\nc\bbN{\mathbb{N}}\nc\bbF{\mathbb{F}}\nc\bbR{\mathbb{R}}
\nc\cA{\mathcal{A}}\nc\cF{\mathcal{F}}\nc\cU{\mathcal{U}}\nc\cI{\mathcal{I}}\nc{\cO}{\mathcal{O}}
\nc\cV{\mathcal{V}}\nc\cN{\mathcal{N}}\nc\cR{\mathcal{R}}
\nc\cS{\mathcal{S}}\nc\cB{\mathcal{B}}\nc\cC{\mathcal{C}}\nc\cM{\mathcal{M}}
\nc\un{\bigcup} \nc\x{\times} \nc\sub{\subseteq} \nc\sps{\supseteq} \nc\sm{\setminus}
\nc{\bset}[2]{\left\{\,#1:#2\,\right\}}\nc{\mset}[2]{\Bigl\{\,#1:#2\,\Bigr\}}
\nc{\set}[2]{\{\,#1:#2\,\}}
\nc{\card}[1]{\left |#1\right |}
\nc{\op}[1]{\operatorname{#1}}
\nc\stab{\op{stab}} \nc\FP{\op{FP}} \nc\FS{\op{FS}}
\nc\dimn{\op{dim}} \nc\spann{\op{span}}\nc{\Fin}{\op{Fin}(\bbN)}\nc{\Finsets}{\op{Fin}}
\nc\AP{\mathrm{AP}}
\nc{\comp}{^{\,\texttt{c}}}
\nc{\cl}[1]{\overline{#1}}
\nc{\Hint}[1]{\par\noindent\emph{Hint}: #1}
\nc{\NN}{\bbN^\bbN}

\newtheorem{thm}{Theorem}[section]\nc{\bthm}{\begin{thm}} \nc{\ethm}{\end{thm}}
\newtheorem{prop}[thm]{Proposition}\nc{\bprp}{\begin{prop}} \nc{\eprp}{\end{prop}}
\newtheorem{fact}[thm]{Fact}\nc{\bfct}{\begin{fact}} \nc{\efct}{\end{fact}}
\newtheorem{prob}[thm]{Problem}\nc{\bprb}{\begin{prob}} \nc{\eprb}{\end{prob}}
\newtheorem{lem}[thm]{Lemma}\nc{\blem}{\begin{lem}} \nc{\elem}{\end{lem}}
\newtheorem{claim}[thm]{Claim}\nc{\bclm}{\begin{claim}} \nc{\eclm}{\end{claim}}
\newtheorem{cor}[thm]{Corollary}\nc{\bcor}{\begin{cor}} \nc{\ecor}{\end{cor}}
\newtheorem{conj}[thm]{Conjecture}\nc{\bcnj}{\begin{conj}} \nc{\ecnj}{\end{conj}}
\theoremstyle{definition}\newtheorem{defn}[thm]{Definition}
\nc{\bdfn}{\begin{defn}} \nc{\edfn}{\end{defn}}
\theoremstyle{remark}
\newtheorem{rem}[thm]{Remark}\nc{\brem}{\begin{rem}}\nc{\erem}{\end{rem}}
\newtheorem{cnv}[thm]{Convention}\nc{\bcnv}{\begin{cnv}} \nc{\ecnv}{\end{cnv}}
\newtheorem{exam}[thm]{Example}\nc{\bexm}{\begin{exam}}\nc{\eexm}{\end{exam}}
\newtheorem{exz}[thm]{Exercise}\nc{\bexz}{\begin{exz}}\nc{\eexz}{\end{exz}}
\nc{\bpf}{\begin{proof}}\nc{\epf}{\end{proof}}\nc{\be}{\begin{enumerate}}\nc{\ee}{\end{enumerate}}
\nc{\bi}{\begin{itemize}}\nc{\itm}{\item}\nc{\ei}{\end{itemize}}
\nc{\beqn}{\begin{eqnarray*}}\nc{\eeqn}{\end{eqnarray*}}

\nc{\ed}{

\subsection*{Acknowledgments}
Marion Scheepers was the first to realize the connection between Ramsey theory and
selection principles, by proving the following beautiful
qualitative extension of Ramsey's Theorem~\cite{coc1}:
Let $X$ be a topological space. If $\sone(\Om(X),\Om(X))$ holds, then for each
cover $\cU\in\Om(X)$ and each finite coloring of the set $[\cU]^2$, there is in $\Om(X)$ a cover $\cV\sub\cU$
such that the graph $[\cV]^2$ is monochromatic. Scheepers proved a large number of results of
this type, including ones jointly with Ljubi\v{s}a Ko\v{c}inac and others (e.g., \cite{coc7, SchForcing, GlCovs}).

Regarding the earlier paper~\cite{suf}, Terence Tao asked me whether
superfilters, viewed as subsets of $\beta S$, are closed.
A positive answer was known, but Tao's question pointed in a fruitful direction.
David J. Fern\'andez Bret\'on, Gili Golan and Michael (Micha\l{}) Machura read
drafts of this paper and made excellent comments.
I had long and helpful correspondences with Neil Hindman and
Imre Leader on the Milliken--Taylor Theorem.
A substantial part of this research was carried out during my Sabbatical leave at
the Faculty of Mathematics and Computer Science, Weizmann Institute of Science. I
thank Gideon Schechtman and the Faculty for their kind hospitality,
and the Faculty Teaching Committee chair, Itai Benjamini,
for the opportunity to deliver a course on Ramsey theory.
This course helped me shape the theory presented here.

\end{document}
}

\title{Algebra, selections, and additive Ramsey theory}

\author{Boaz Tsaban}

\address{Department of Mathematics, Bar-Ilan University, Ramat Gan 5290002, Israel,
and
Faculty of Mathematics and Computer Science, Weizmann Institute of Science, Rehovot 76100\-01, Israel.}

\email{tsaban@math.biu.ac.il}

\urladdr{math.biu.ac.il/~tsaban}

\begin{document}

\begin{abstract}
Hindman's celebrated Finite Sums Theorem, and its high-dimensional version due to Milliken and Taylor,
are extended from covers of countable sets
to covers of arbitrary topological spaces with Menger's  classic covering property.
The methods include, in addition to Hurewicz's game theoretic characterization of Menger's
property, extensions of the classic idempotent theory in the Stone--\v{C}ech compactification 
of semigroups, and of the more recent theory of selection principles.
This provides strong versions of the mentioned celebrated theorems, where the monochromatic substructures
are large, beyond infinitude, in an analytic sense.
Reducing the main theorems to the purely combinatorial setting, we obtain nontrivial consequences concerning uncountable cardinal characteristics of the continuum.

The main results, modulo technical refinements, are of the following type
(definitions provided in the main text):
Let $X$ be a Menger space, and $\mathcal{U}$ be an infinite open cover of $X$.
Consider the complete graph, whose vertices are the open sets in $X$.
For each finite coloring of the vertices and edges of this graph,
there are disjoint finite subsets $\mathcal{F}_1,\mathcal{F}_2,\dots$ of the cover $\mathcal{U}$
whose unions $V_1 := \bigcup\mathcal{F}_1, V_2 := \bigcup\mathcal{F}_2,\dots$
have the following properties:
\begin{enumerate}
\item The sets $\bigcup_{n\in F}V_n$ and $\bigcup_{n\in H}V_n$ are distinct for all nonempty finite sets
$F<H$.
\item All vertices $\bigcup_{n\in F}V_n$, for nonempty finite sets $F$, are of the same color.
\item All edges $\bigl\{\,\bigcup_{n\in F}V_n, \bigcup_{n\in H}V_n\,\bigr\}$, for nonempty
finite sets $F<H$, have the same color.
\item The family $\{V_1,V_2,\dots\}$ is an open cover of $X$.
\end{enumerate}
A self-contained introduction to the necessary parts of the needed theories is provided.
\end{abstract}



\maketitle

\section{Background}

Following is a brief, self-contained introduction to the Stone--\v{C}ech compactification of a
semigroup and its necessary algebraic and combinatorial properties. All assertions made can be
verified directly. More detailed introductions, with additional
combinatorial applications, may be found in the books \cite{HS, ProtasovRT}.
Familiar parts may be skipped by the reader.

\subsection{The Stone--\v{C}ech compactification and Hindman's Theorem}
Almost throughout, $S$ denotes an infinite semigroup.
We do not assume that the semigroup $S$ is commutative; however, with the applications in mind,
we use additive notation.
The \emph{Stone--\v{C}ech compactification} of $S$, $\beta S$, is the set of all ultrafilters on $S$.
We identify each element $s\in S$ with the principal ultrafilter associated to it.
Thus, we view the set $S$ as a subset of $\beta S$.
A filter $\cF$ on $S$ is \emph{free} if the intersection $\bigcap\cF$ of all elements of $\cF$ is empty.
An \emph{ultrafilter} is free if and only if it is nonprincipal.

A topology on the set $\beta S$ is defined by taking the sets
$[A]:=\set{p\in\beta S}{A\in p}$, for $A\sub S$, as a basis for the topology.
The function $A\mapsto [A]$ respects finite unions, finite intersections, and complements.
For an element $s\in S$ and a set $A\sub S$, we have that $s\in [A]$ if and only if $s\in A$.
In particular,  the set $S$ is dense in $\beta S$.

The topological space $\beta S$ is compact: If $\beta S=\Un_{\alpha\in I}[A_\alpha]$
and no finite union of sets $A_\alpha$ is $S$, then the family $\set{A_\alpha\comp}{\alpha\in I}$
extends to an ultrafilter $p\in\beta S$, so $p$ is in some set $[A_\alpha]$; a contradiction.
Define the sum of elements $p,q\in\beta S$ by
$$A\in p+q\mbox{ if and only if }\set{b\in S}{\exists C\in q,\ b+C\sub A}\in p.$$
Then $p+q\in\beta S$. We obtain an extension of the addition operator from $S$ to $\beta S$,
with the following continuity properties:
\be
\item For each element $x\in S$,
the function $q\mapsto x+q$ is continuous.
\item For each element $q\in \beta S$, the function $p\mapsto p+q$ is continuous.
\ee
Fix $x,y\in S$. Since $(x+y)+z=x+(y+z)$ for all $z\in S$ and the set $S$ is dense is $\beta S$,
we have by (1) that $(x+y)+r=x+(y+r)$ for all $r\in\beta S$.
Fixing $r$ and unfixing $y$, we have by (1) and (2) that $(x+q)+r=x+(q+r)$
for all $q\in\beta S$. Finally, fixing $q$ and unfixing $x$, we have by (2)
that $(p+q)+r=p+(q+r)$ for all $p\in\beta S$. Thus, $(\beta S,+)$ is a semigroup.

If $e\in\beta S$ is an idempotent element, that
is, if $e+e=e$, then for each set $A\in e$ there are a set $B\in e$, and for each $b\in B$,
a set $C\in e$ such that $b+C\sub A$.
Conversely, the latter property of $e$ implies that $e\sub e+e$ and thus $e=e+e$.
In this characterization, by intersecting $C$ with $A$, we may assume that $C\sub A$.

By the continuity of the functions $p\mapsto p+q$, for $q\in\beta S$,
there are idempotent elements in any closed subsemigroup $T$ of $\beta S$.
Indeed, Zorn's Lemma
provides us with a minimal closed subsemigroup $E$ of $T$, and it follows by minimality
that $E=\{e\}$ for some (necessarily, idempotent) element $e\in T$.\footnote{To see that a minimal
closed subsemigroup $E$ of $\beta S$ must be of the form $\{e\}$,
fix an element $e\in E$. As $E+e$ is a closed subsemigroup of $E$, we have that $E+e=E$.
Thus, the stabilizer of $e$, $\set{t\in E}{t+e=e}$ is a (closed) subsemigroup of $E$, and is therefore equal to $E$. Then $e+e=e$.
}

\bdfn
For elements $\seq{a}$ in a semigroup $S$,
and a nonempty finite set $F=\{i_1,\dots,i_k\}\sub\bbN$ with $k\ge 1$ and
$i_1<\cdots<i_k$, define
$a_F:=a_{i_1}+\cdots+a_{i_k}$.
Let
$$\FS(\seq{a}):=\set{a_F}{F\sub\bbN,\mbox{ $F$ finite nonempty}},$$
the set of all finite sums, in increasing order of indices, of elements $a_i$.
Similarly, for elements $a_1,\dots,a_n\in S$, the set $\FS(a_1,\dots,a_n)$ is
comprised of the elements $a_F$ for $F$ a nonempty
subset of $\{1,\dots,n\}$.
\edfn

A \emph{finite coloring} of a set $A$ is a function $f\colon A\to\{1,\dots,k\}$, for $k\in\bbN$.
Given a finite coloring $f$ of a set $A$, a set $B\sub A$ is \emph{monochromatic} if
there is a color $i$ with $f(b)=i$ for all $b\in B$.

\bthm[Hindman \cite{Hindman74}]
For each finite coloring of $\bbN$, there are elements $\seq{a}\in\bbN$ such that the set
$\FS(\seq{a})$ is monochromatic.
\ethm

The following strikingly elegant proof of Hindman's Theorem is due to Galvin and Glazer.
Fix an idempotent element $e\in\beta\bbN$.
Let a $k$-coloring of $\bbN$ be given.
If $C_i$ is the set of elements of color $i$, then $C_1\cup\cdots\cup C_k=\bbN\in e$, and thus there is
a color $i$ with $A_1:=C_i\in e$.
For $n=1,2,\dots$, since $e$ is an idempotent ultrafilter,
there are an element $a_n\in A_n$ and a set $A_{n+1}\sub A_n$ in $e$
such that $a_n+A_{n+1}\sub A_n$. It then follows, considering the sums from right to left, that
every finite sum $a_{i_1}+\cdots+a_{i_k}$, for $i_1<\cdots<i_k$, is in $A_{i_1}$.
Thus, the set $\FS(\seq{a})$ is a subset of the monochromatic set $A_1$.

\subsection{The Milliken--Taylor Theorem and proper sumsequences}

For a set $S$, let $[S]^2$ be the set of all $2$-element subsets
of $S$; equivalently, the edge set of the complete graph with vertex set $S$.

\bdfn
Let $S$ be a semigroup.
For nonempty finite sets of natural numbers $F$ and $H$, we write
$F<H$ if all elements of $F$ are smaller than all elements of $H$.
A \emph{sumsequence} (or sum subsystem)
of a sequence $\seq{a}\in S$ is a sequence of the form
$a_{F_1},a_{F_2},\dots$, for nonempty finite sets of natural numbers $F_1<F_2<\cdots$.

A sequence $\seq{b}\in S$ is \emph{proper} if $b_F\neq b_H$
for all nonempty finite sets $F<H$ of natural numbers.
The \emph{sum graph} of a proper sequence $\seq{b}\in S$
is the subset of $[\FS(\seq{b})]^2$ consisting of the
edges $\{b_F, b_H\}$, for nonempty finite sets $F<H$ of natural numbers.
\edfn

If $\seq{b}$ is a sumsequence of $\seq{a}$, then $\FS(\seq{b})\sub\FS(\seq{a})$.
The relation of being a sumsequence is transitive.

Ramsey's Theorem~\cite{Ramsey} asserts that, for each finite coloring of an infinite complete graph $[V]^2$
with vertex set $V$, there is an infinite complete monochromatic subgraph, that is, an infinite set
$I\sub V$ such that the set $[I]^2$ is monochromatic.
The Milliken--Taylor Theorem unifies Hindman's and Ramsey's theorems.

\bthm[Milliken--Taylor \cite{Milliken, Taylor}]\label{thm:MTN}
Let $\seq{a}$ be a sequence in $\bbN$.
For each finite coloring of the set $[\bbN]^2$, there is a proper
sumsequence $\seq{b}$ of the sequence $\seq{a}$ such that
the sum graph of $\seq{b}$ is monochromatic.
\ethm

The Milliken--Taylor Theorem can be proved by combining the
proofs of Ramsey's and Hindman's Theorems, as
can be gleaned from the proof of the forthcoming Theorem~\ref{thm:abs}.

Our applications are in a setting where all elements of the semigroup $S$ are idempotents.
In this case, stating Hindman's Theorem for the semigroup $S$ instead of $\bbN$ yields a trivial statement:
for an idempotent element $e\in S$,
the set $\FS(e,e,\dots)=\{e\}$ is obviously monochromatic.
The sequence $e,e,\dots$ is improper, and so are all of its sumsequences.
Thus, the Milliken--Taylor Theorem cannot be extended to such cases.
An example of a semigroup with all elements idempotent is $\Fin$, the set of nonempty
finite subsets of $\bbN$, with the operation $\cup$. For this semigroup, we have the following.

\bthm[Milliken--Taylor]\label{thm:MTFin}
For each finite coloring of the set $[\Fin]^2$, there
are elements $F_1<F_2<\cdots$ in $\Fin$ such that
the sum graph of $\seq{F}$ is monochromatic.
\ethm

As every sequence of natural numbers has a proper sumsequence,
Theorem~\ref{thm:MTN} is a special case of the following one.

\bthm\label{thm:MTS}
Let $S$ be a semigroup, and $\seq{a}\in S$. If the sequence $\seq{a}$
has a proper sumsequence, then for each finite coloring of the set $[S]^2$, there is
a proper sumsequence of $\seq{a}$ whose sum graph is monochromatic.
\ethm

Theorem~\ref{thm:MTS} is more general than Theorem~\ref{thm:MTFin},
since the sequence $\{1\},\{2\},\allowbreak\dots$ is proper.
Theorem~\ref{thm:MTS} follows from Theorem~\ref{thm:MTFin}:
By moving to a sumsequence, we may assume that the sequence $\seq{a}$ is proper.
Let $\chi$ be a finite coloring of $[S]^2$. Define a coloring $\kappa$ of $[\Fin]^2$ by
$\kappa(\{F,H\})=\chi(\{a_F,a_H\})$ for $F<H$, and $\kappa(\{F,H\})$ arbitrary otherwise,
and apply Theorem~\ref{thm:MTFin},
using that sumsequences of proper sequences are proper.

The hypothesis of having a proper sumsequence fails only in degenerate cases.

\bprp\label{prp:proper}
Let $S$ be a semigroup, and $\seq{a}\in S$. If the sequence $\seq{a}$ has no
proper sumsequence, then every sumsequence of $\seq{a}$ has a sumsequence
of the form $e,e,\dots$, where $e$ is an idempotent element of $S$; equivalently,
a sumsequence whose set of finite sums is a singleton.
\eprp
\bpf
We use Theorem~\ref{thm:MTFin}.
Define a coloring of the set $[\Fin]^2$ by
$$\{F,H\}\mapsto\card{\{a_F,a_H\}}.$$
Let $F_1<F_2<\cdots$ be elements of $\Fin$
such that the sum graph of the sequence $\seq{F}$ is monochromatic.

Consider the sumsequence $b_1:=a_{F_1},b_2:=a_{F_2},\dots$.
Assume that the color is 2. Then the sumsequence $\seq{b}$ is proper; a contradiction.
Thus, the color must be 1.
Then $b_{H_1}=b_{H_2}$ for all $H_1<H_2$ in $\Fin$.
Let $e:=b_1$. Then $b_n=b_1=e$ for all $n>1$. For each set $H\in\Fin$, take $n>H$.
Then $b_H=b_n=e$. In particular, $e+e=b_1+b_2=b_{\{1,2\}}=e$.
Thus, $\FS(\seq{b})=\{e\}$.
\epf

\section{Idempotent filters and superfilters}

Superfilters provide a convenient way to identify closed subsets of $\beta S$.\footnote{Superfilters
have various names in the classic literature, including
coideals, grilles, and partition-regular families, depending on the context where they are used.
Some of the definitions in the literature are not equivalent to the one given here, but they are
always conceptually similar. The present term is adopted from~\cite{suf}.}

\bdfn
A family $\cA$ of subsets of a set $S$ is a \emph{superfilter} on $S$ if:
\be
\item All sets in $\cA$ are infinite.
\item For each set $A\in\cA$, all subsets of $S$ that contain $A$ are in $\cA$.
\item Whenever $A_1\cup A_2\in\cA$, we have that $A_1$ or $A_2$ are in $\cA$;
equivalently, for each set $A\in\cA$ and each finite coloring of $A$,
there is in $\cA$ a monochromatic subset of $A$.
\ee
\edfn

The simplest example of a superfilter on a set $S$ is the family $[S]^\infty$,
consisting of all infinite subsets of $S$.
Many examples of superfilters are provided by Ramsey theoretic theorems. For example,
van der Waerden's Theorem asserts that monochromatic arithmetic progressions of any prescribed
finite length will be found in any long enough, finitely-colored arithmetic progression.
By van der Waerden's Theorem, the family of all sets of natural numbers
containing arbitrarily long finite arithmetic progressions is a superfilter on $\bbN$.

The notions of free filter and superfilter are dual.
For a family $\cF$ of subsets of a set $S$, define
$\cF^+:=\set{A\sub S}{A\comp\notin\cF}$.
The following assertions are easy to verify.

\blem[Folklore]
\label{lem:folk}
Let $S$ be a set.
\be
\item For all families $\cF_1$ and $\cF_2$ of subsets of $S$, $\cF_1\sub\cF_2$ implies
that $\cF_1^+\sps\cF_2^+$.
\item For each family $\cF$ of subsets of $S$, $\cF^{++}=\cF$.
\item For each free filter $\cF$ on $S$, the set $\cF^+$ is a superfilter containing $\cF$.
\item For each superfilter $\cA$ on $S$, the set $\cA^+$ is a free filter contained in $\cA$.
\item For each filter $\cF$, if $A\in\cF^+$ and $B\in\cF$, then $A\cap B\in\cF^+$.
\item For each ultrafilter $p$ on $S$, $p^+=p$.
\ee
\elem
\bpf[Proof of (5)]
Since $A\sub B\comp\cup (A\cap B)$, the latter set is in $\cF^+$.
Since $B\comp\notin\cF^+$, we have that $A\cap B\in\cF^+$.
\epf

Every free ultrafilter on $S$ is a superfilter on $S$, and so is any union
of free ultrafilters on $S$.
Since elements of superfilters are infinite, the filter of cofinite subsets of $S$
is contained in all superfilters on $S$.
By the following lemma, every superfilter $\cA$ is a union of a closed set of free ultrafilters.
Indeed, taking $\cF=\{\bbN\}$ we have by the lemma that the set
$C:=\set{p\in\beta S}{p\sub\cA}$ is closed,
and for each set $A\in\cA$, letting $\cF$ be the filter generated by $A$ we see, again by the lemma,
that there is an ultrafilter $p\in C$ with $A\in p$. Thus, $\Un C=\cA$.

\blem\label{lem:Fsuf}
Let $S$ be an infinite set.
For each superfilter $\cA$ on $S$, and each filter $\cF\sub\cA$,
the set $\set{p\in\beta S}{\cF\sub p\sub\cA}$ is a nonempty closed
subset of $\beta S\sm S$.
\elem
\bpf
It is straightforward to verify that the set is closed. We prove that it is nonempty.
By Lemma~\ref{lem:folk}(4), the set $\cA^+$ is a filter.
By Lemma~\ref{lem:folk}(2,5) applied to the filter $\cA^+$,
we have that $A\cap B\in\cA$ for all $A\in\cA$ and $B\in\cA^+$.
In particular, the set $A\cap B$ is infinite for all $A\in\cF, B\in\cA^+$.
The family $\set{A\cap B}{A\in\cF, B\in\cA^+}$ is closed under finite intersections.
Since its elements are infinite, it extends to a free ultrafilter $p$.
Necessarily, $\cF\sub p$. If there were an element $B\in p\sm \cA$, then
$B\comp\in\cA^+\sub p$; a contradiction.
\epf

\bdfn
Let $S$ be a semigroup.
\be
\item For a set $A\sub S$ and a family $\cF$ of subsets of $S$, let
$$A^\star(\cF):=\set{b\in S}{\exists C\in\cF, b+C\sub A}.$$
\item A filter $\cF$ on $S$ is an \emph{idempotent filter}
if for each set $A\in\cF$, the set $A^\star(\cF)$ is in $\cF$.
\item A superfilter $\cA$ on $S$ is an \emph{idempotent superfilter}
if, for each set $A\sub S$ such that the set $A^\star(\cA)$ is in $\cA$,
we have that  $A\in\cA$.
\ee
\edfn

Thus, for ultrafilters $p,q$ on $S$, $A\in p+q$ if and only if $A^\star(q)\in p$.

Let $S$ be a semigroup.
A superfilter $\cA$ on $S$ is \emph{translation-invariant}
if $s+A\in\cA$ for all $s\in S$ and $A\in\cA$.
Every translation-invariant superfilter on a semigroup $S$ is an idempotent superfilter.

Since ultrafilters are maximal filters, we have that, for an ultrafilter $p$ on a semigroup $S$,
being an idempotent ultrafilter, idempotent filter, and idempotent superfilter is the same.

\blem\label{lem:FAplus}
Let $S$ be a semigroup.
\be
\item For each free idempotent filter $\cF$ on $S$, the superfilter $\cF^+$ is idempotent.
\item For each idempotent superfilter $\cA$ on $S$, the free filter $\cA^+$ is idempotent.
\ee
\elem
\bpf
(1) Let $A\sub S$, and assume that the set $B_1:=A^\star(\cF^+)$ is in $\cF^+$.
Assume that $A\notin\cF^+$.
Then $A\comp\in\cF$, and thus the set $B_2:=(A\comp)^\star(\cF)$ is in $\cF$.
By Lemma~\ref{lem:folk}(5), there is an element $b\in B_1\cap B_2$.
Then there are sets $C_1\in\cF^+$ and $C_2\in\cF$ such that $b+C_1\sub A$
and $b+C_2\sub A\comp$. Pick $c\in C_1\cap C_2$. Then $b+c\in A\cap A\comp$; a contradiction.

(2) Similar.
\epf

\blem\label{lem:ipextend}
Let $S$ be a semigroup, and $\cF$ be a  free idempotent filter on $S$.
Then the set $T:=\set{p\in\beta S}{\cF\sub p}$ is a closed subsemigroup of $\beta S$ disjoint from $S$.
\elem
\bpf
By Lemma~\ref{lem:Fsuf}, with $\cA=[S]^\infty$, the set $T$ is a closed subset of $\beta S$.
Since the filter $\cF$ is free, we have that $T\sub\beta S\sm S$.
Let $p,q\in T$, and $A\in\cF$.
Since the filter $\cF$ is idempotent, $A^\star(\cF)\in\cF\sub p$.
Since $\cF\sub q$, $A^\star(\cF)\sub A^\star(q)$, and therefore $A^\star(q)\in p$.
By the definition of sum of ultrafilters, $A\in p+q$.
\epf

\bthm\label{thm:ip}
Let $S$ be a semigroup, and assume that $\cF$ is a free idempotent filter on $S$
contained in an idempotent superfilter $\cA$ on $S$.
Then there is a free idempotent ultrafilter $e$ with $\cF\sub e\sub\cA$.
\ethm
\bpf
Let $T_1=\set{p\in\beta S}{\cF\sub p}$ and $T_2:=\set{p\in\beta S}{p\sub\cA}$.
By Lemma~\ref{lem:ipextend}, the set $T_1$ is a closed subsemigroup of $\beta S$,
and so is the set $\set{p\in\beta S}{\cA^+\sub p^+=p}=T_2$.

By Lemma~\ref{lem:Fsuf}, the intersection $T:=T_1\cap T_2$ is nonempty, and is therefore a closed
subsemigroup of $\beta S$. Pick an idempotent element in $T$.
\epf

\section{Selection principles and an abstract partition theorem}

We use the following notions from Scheepers's seminal paper~\cite{coc1}.
Let $\cA$ and $\cB$ be families of sets.
$\sone(\cA,\cB)$ is the property that, for each sequence $\seq{A}\in\cA$,
one can select one element from each set,
$b_1\in A_1,b_2\in A_2,\dots$,
such that $\{\seq{b}\}\in\cB$.
$\gone(\cA,\cB)$ is a game associated to $\sone(\cA,\cB)$.
This game is played by two players, \emph{Alice} and \emph{Bob},
and has an inning per each natural number.
In the $n$-th inning, Alice plays a set $A_n\in\cA$, and Bob selects an element $b_n\in A_n$.
Bob wins if $\{\seq{b}\}\in\cB$. Otherwise, Alice wins.

If Alice does not have a winning strategy in the game $\gone(\cA,\cB)$, then
$\sone(\cA,\cB)$ holds. The converse implication holds in some important cases,
including the ones in our main applications.
A survey of known results of this type is provided, e.g., in Section~11 of~\cite{LecceSurvey}.

\bexm\label{exm:FplusGame}
Let $S$ be a set, and $\cF$ be a filter on $S$ generated by countably many sets.
Then Alice does not have a winning strategy in the game $\gone(\cF^+,\cF^+)$;
moreover, Bob has one: Fix sets $\seq{B}\in\cF$ such that every member
of $\cF$ contains one of these sets.
In each inning, by Lemma~\ref{lem:folk}(5),
Bob can pick an element $b_n\in A_n\cap B_n$. Then $\{\seq{b}\}\in\cF^+$.

For the filter $\cF$ of cofinite sets, this reproduces the simple observation
that Bob has a winning strategy in the game $\gone([S]^\infty,[S]^\infty)$.
\eexm

In general, the game $\gone(\cA,\cB)$ is not determined,
and the property that Alice does not have a winning strategy is strictly weaker than Bob's having one.
This will be the case in our main applications \cite[Section~11]{LecceSurvey}.

\bdfn
A \emph{free idempotent chain} in a semigroup $S$ is a descending sequence
$A_1\sps A_2\sps \cdots$ of infinite subsets of $S$ such that:
\be
\item $\bigcap_n A_n=\emptyset$.
\item For each $n$, the set $A_n^\star(\{\seq{A}\})$
contains one of the sets $A_m$; equivalently, there is $m>n$ such that,
for each $a\in A_m$, there is $k>m$ with $a+A_k\sub A_m$.
\ee
For a family $\cA$ of subsets of $S$, a \emph{free idempotent chain in $\cA$}
is a free idempotent chain of elements of $\cA$.
\edfn

\bexm\label{exm:chain}
For each proper sequence $\seq{a}$ in a semigroup, the sets
$\FS(a_n,a_{n+1},\dots)$, for $n\in\bbN$, form a free idempotent chain.
Thus, if a sequence $\seq{a}$ has a \emph{proper sumsequence},
then there is a free idempotent chain $A_1\sps A_2\sps\cdots$
with $A_n\sub\FS(a_n,a_{n+1},\dots)$ for all $n$.
\eexm

\blem\label{lem:chainfilter}
Let $S$ be a semigroup, and $\cA$ be a superfilter on $S$.
Every filter generated by a free idempotent chain in $\cA$ is a free idempotent filter contained in $\cA$.
\elem
\bpf
Let $A_1\sps A_2\sps \cdots$ be a free idempotent chain in $\cA$, and let $\cF$ be the filter
generated by the sets $\seq{A}$.
Since $\bigcap_nA_n=\emptyset$, the filter $\cF$ is free.
Since $A_n\in\cA$ for each $n$, $\cF\sub\cA$.
The filter $\cF$ is idempotent: For $A\in\cF$, let $A_n$ be a subset of $A$.
By the definition, there is $m$ such that $A^\star(\cF)\sps A_n^\star(\{\seq{A}\})\sps A_m$.
Since $A_m\in\cF$, we have that $A^\star(\cF)\in\cF$.
\epf

Our theorems can be stated for any finite dimension. For clarity,
we state them in the one-dimensional case, which extends Hindman's Theorem,
and in the two-dimensional case, which extends the Milliken--Taylor Theorem.
The one-dimensional case always follows from the two-dimensional, for the following
reason.

\bprp\label{prp:2dim1dim}
Let $S$ be a semigroup, and $\chi$ be a finite coloring of the sets $S$ and $[S]^2$.
There is a finite coloring $\eta$ of the set $[S]^2$ such that,
for each proper sequence $\seq{b}$ with $\eta$-monochromatic sum graph,
the set $\FS(\seq{b})$ and the sum graph of $\seq{b}$ are both $\chi$-monochromatic.
\eprp
\bpf
By enumerating the elements of the countable set $\FS(\seq{b})$,
we obtain an order on this set such that every element has
only finitely many smaller elements.
Define a coloring $\kappa$ of $[\FS(\seq{b})]^2$ by
$$\kappa(\{s,t\}):=\chi(\min\{s,t\}).$$
Extend $\kappa$ to a coloring of $[S]^2$ in an arbitrary manner.

Assume that the set $\FS(\seq{b})$ is monochromatic for $\kappa$, say green.
Being proper, the sequence $\seq{b}$ is bijective.
For each nonempty finite set $F$ of natural numbers,
since there are at most finitely many elements in $\FS(\seq{b})$ smaller
than $b_F$, there is $n>F$ such that $b_F<b_n$. Then $\kappa(\{b_F,b_n\})=\chi(b_F)$.
Thus, the element $b_F$ is green.

The finite coloring $\eta$ of the set $[S]^2$, defined by
$$\eta(\{s,t\}):=\bigl(\,\kappa(\{s,t\}),\chi(\{s,t\})\,\bigr),$$
is as required.
If $\chi$ is a $k$-coloring, we may represent the range set of $\eta$ in the form $\{1,\dots,k^2\}$.
\epf

The two monochromatic sets in Proposition~\ref{prp:2dim1dim} may be of different colors.
Moreover, this can be forced by adding a coordinate to $\chi(x)$ that is
$1$ if $x\in S$ and $2$ if $x\in[S]^2$.

\bthm\label{thm:abs}
Let $S$ be a semigroup.
Let $\cA$ be an idempotent superfilter on $S$, and $\cB$ be a family of subsets of $S$
such that Alice does not have a winning strategy in the game $\gone(\cA,\cB)$.
Let $\seq{a}$ be a sequence in $S$, and $A_1\sps A_2\sps\cdots$ be a free
idempotent chain in $\cA$ with $A_n\sub\FS(a_n,a_{n+1},\dots)$ for all $n$.
For each finite coloring of the sets $S$ and $[S]^2$,
there are elements $b_1\in A_1, b_2\in A_2,\dots$ such that:
\be
\item The set $\{\seq{b}\}$ is in $\cB$.
\item The sequence $\seq{b}$ is a proper sumsequence of $\seq{a}$.
\item The set $\FS(\seq{b})$ is monochromatic.
\item The sum graph of $\seq{b}$ is monochromatic.
\ee
\ethm
\bpf
By Proposition~\ref{prp:2dim1dim}, it suffices to prove the two-dimensional assertion, that is,
item (3) follows from item (4).

By Lemma~\ref{lem:chainfilter}, there is a free idempotent filter $\cF$ such that
$\{\seq{A}\}\sub\cF\sub\cA$.
By Theorem~\ref{thm:ip}, there is a free idempotent ultrafilter $e$ on $S$ such that
$\cF\sub e\sub \cA$.

Let a finite coloring $\chi\colon [S]^2\to \{1,\dots,k\}$ be given.
For each element $s\in S$, let
$$C_i(s):=\set{t\in S\sm\{s\}}{\chi(\{s,t\})=i}.$$
As $C_1(s)\cup\cdots C_k(s)=S\sm\{s\}\in e$, there is a unique $i$ with $C_i(s)\in e$.
Define a finite coloring $\kappa\colon S\to\{1,\dots,k\}$ by letting $\kappa(s)$ be
this unique $i$ with $C_i(s)\in e$.
Since $e$ is an ultrafilter, there is in $e$ a set $M\sub S$
that is monochromatic for the coloring $\kappa$.
Assume that the color is \emph{green}.
Then, for each finite set $F\sub M$, we have that
$$G(F):=\bigcap_{s\in F}\set{t\in S\sm\{s\}}{\{s,t\}\mbox{ is green}}\in e,$$
and for each element $s\in F$ and each element $t\in G(F)$,
we have that $s\neq t$ and the edge $\{s,t\}$ is green.

For a set $D\in e$, define
$$D^\star:=\set{b\in D}{\exists B\sub D\mbox{ in }e, b+B\sub D}=D^\star(e)\cap D.$$
Then $D^\star\sub D$ and, since $e$ is an idempotent ultrafilter, $D^\star\in e$.

We define a strategy for Alice. In this strategy, Alice makes choices from certain nonempty sets.
Formally, she does that by applying prescribed choice functions to the given nonempty sets.
\be
\item In the first inning, Alice sets $D_1:=M\cap A_1$,
and plays the set $D_1^\star$.

\item Assume that Bob plays an element $b_1\in D_1^\star$.
Then Alice chooses a set $B\sub D_1$ in $e$ such that $b_1+B\sub D_1$
and a set $F_1$ with $a_{F_1}=b_1$.
She then chooses a natural number $m_1>F_1$, and sets
$D_2:=B\cap G(\{b_1\})\cap A_{m_1}$.
Having done that, Alice plays the set $D_2^\star$.

\item Assume that Bob plays an element $b_2\in D_2^\star$.
Then $b_1+b_2\in D_1\sub M$.
Alice chooses a set $B\sub D_2$ in $e$ such that $b_2+B\sub D_2$,
a set $F_2>m_1$ with $a_{F_2}=b_2$,
and a natural number $m_2>F_2$.
She sets $D_3:=B\cap G(\FS(b_1,b_2))\cap A_{m_2}$, and plays $D_3^\star$.

\item In the $n+1$-st inning, Bob has picked elements $b_1\in D_1^\star,\dots,b_n\in D_n^\star$.
As in the Galvin--Glazer proof of Hindman's Theorem,
by computing sums from right to left, we see that $\FS(b_1,\dots,b_n)\sub D_1\sub M$.
Alice chooses a set $B\sub D_n$ in $e$ such that $b_n+B\sub D_n$,
a set $F_n>m_{n-1}$ with $a_{F_n}=b_n$,
and a natural number $m_n>F_n$.
She then sets $D_{n+1}:=B\cap G(\FS(b_1,\dots,b_n))\cap A_{m_n}$,
and plays the set $D_{n+1}^\star$.
\ee
Since Alice has no winning strategy, there is a play $(D_1^\star,b_1,D_2^\star,b_2,\dots)$,
according to Alice's strategy, won by Bob.
By the construction, the sequence $\seq{b}$ is a sumsequence of $\seq{a}$.
The set $\{\seq{b}\}$ is in $\cB$, since Bob won this play.

Let $i_1<\cdots<i_k<j_1<\cdots<j_l$, $F=\{i_1,\dots,i_k\}$, and $H=\{j_1,\dots,j_l\}$.
Then
\begin{eqnarray*}
b_{F} & \in & \FS(b_1,\dots,b_{i_k}), \mbox{ and}\\
b_{H} & = & b_{j_1}+\cdots+b_{j_l}.
\end{eqnarray*}
Computing the latter sum from right to left, we see that
$$b_{H}\in D_{j_1}\sub D_{i_k+1}\sub G(\FS(b_1,\dots,b_{i_k})).$$
It follows that the elements $b_F$ and $b_H$ are distinct,
and the edge $\{b_{F}, b_{H}\}$ is green.
\epf

To gain some intuition on Theorem~\ref{thm:abs}, we provide several simple examples.
These examples can also be established via somewhat more direct arguments.
Theorem~\ref{thm:abs} will be an important part in the proofs of our later, main theorems.

\bexm\label{exm:MT++}
Let $S$ be a semigroup.
Let $\seq{a}$ be a sequence in $S$, and $A_1\sps A_2\sps\cdots$ be a free
idempotent chain with $A_n\sub\FS(a_n,a_{n+1},\dots)$ for all $n$.
For each finite coloring of the sets $S$ and $[S]^2$,
there are elements $b_1\in A_1, b_2\in A_2, \dots$ such that:
\be
\item The sequence $\seq{b}$ is a proper sumsequence of $\seq{a}$.
\item The set $\FS(\seq{b})$ is monochromatic.
\item The sum graph of $\seq{b}$ is monochromatic.
\ee
\eexm
\bpf
By Lemma~\ref{lem:chainfilter}, with the trivial superfilter $\cA=[S]^\infty$,
the filter $\cF$ on $S$ generated by the sets  $\seq{A}$ is a free idempotent filter.
By Lemma~\ref{lem:FAplus}, the superfilter $\cF^+$ is also idempotent.
By Example~\ref{exm:FplusGame}, Bob has a winning strategy in the game $\gone(\cF^+,\cF^+)$.
Since $\cF\sub\cF^+$, Theorem~\ref{thm:abs} applies with $\cA=\cB=\cF^+$.
\epf

In most semigroups $S$ one encounters, left addition is at most finite-to-one. In this case,
the superfilter $[S]^\infty$ is translation-invariant; in particular, idempotent.
In this case, the proof of Example~\ref{exm:MT++} reduces to one short sentence:
Apply Theorem~\ref{thm:abs} with $\cA=\cB=[S]^\infty$.


The Milliken--Taylor Theorem in arbitrary semigroups
(Theorem~\ref{thm:MTS}) follows from Example~\ref{exm:MT++},
by Example~\ref{exm:chain}.

\bexm\label{exm:absQMT}
Let $\seq{\cF}\sub\Fin$, and $A\sub\bbN$.
Assume that every cofinite subset of $A$ contains a member from each family $\cF_n$.
For each finite coloring of the sets $\Fin$ and $[\Fin]^2$,
there are nonempty finite subsets $F_1<F_2<\cdots$ of $A$
such that:
\be
\item Each set $F_n$ contains some element of the family $\cF_n$.
\item All nonempty finite unions $H$ of sets $F_n$ have the same color.
\item All sets $\{H_1,H_2\}$, for $H_1<H_2$ nonempty finite unions of sets $F_n$,
have the same color.
\ee
\eexm
\bpf
We work with the semigroup $\Finsets(A)$ of all nonempty finite subsets of $A$.
Enumerate $A=\{\seq{a}\}$.
For each $n$, let
$$A_n := \set{F\in \Finsets(\{a_n,a_{n+1},\dots\})}{\exists H\in \cF_n, H\sub F}\sub\FS(\{a_n\},\{a_{n+1}\},\dots).$$
Then $\bigcap_n A_n=\emptyset$.
Every set $A_n$ is a subsemigroup of $S$. Thus, the sequence $A_1\sps A_2\sps\cdots$ is a free idempotent
chain.
Apply Example~\ref{exm:MT++}.
\epf


\bexm
Let $A\sub\bbN$ be a set containing arbitrarily long arithmetic progressions.
For each finite coloring of the sets $\Fin$ and $[\Fin]^2$,
there are nonempty finite subsets $F_1<F_2<\cdots$ of $A$ such that:
\be
\item The set $\Un_n F_n$ contains arbitrarily long arithmetic progressions.
\item All nonempty finite unions $H$ of sets $F_n$ have the same color.
\item All sets $\{H_1,H_2\}$, for $H_1<H_2$ nonempty finite unions of sets $F_n$,
have the same color.
\ee
\eexm

Additional examples are provided by any notion that is captured
by finite sets, e.g., entries of solutions of homogeneous systems of equations,
and entries of image vectors of matrices.
The \emph{upper density} of a set $A\sub\bbN$ is the real number
$\limsup_n\card{A\cap\{1,\dots,n\}}/n$.

\bexm
Let $A\sub\bbN$ be a set of upper density $\delta$.
For each finite coloring of the sets $\Fin$ and $[\Fin]^2$,
there are nonempty finite subsets $F_1<F_2<\cdots$ of $A$ such that:
\be
\item The set $\Un_n F_n$ has upper density $\delta$.
\item All nonempty finite unions $H$ of sets $F_n$ have the same color.
\item All sets $\{H_1,H_2\}$, for $H_1<H_2$ nonempty finite unions of sets $F_n$,
have the same color.
\ee
\eexm
\bpf
The upper density of a set does not change by removing finitely many elements from that set.
Take a sequence $\seq{\delta}$ increasing to $\delta$.
For each $n$, let $\cF_n:=\set{F\in\Finsets(A)}{|F|/\max F>\delta_n}$.
Apply Example~\ref{exm:absQMT}.
\epf

An analogous assertion also holds for the so-called \emph{Banach density}.

\section{Menger spaces}

A topological space $X$ is a \emph{Menger space} if, for each sequence $\seq{\cU}$ of
open covers of $X$, there are finite sets $\cF_1\sub\cU_1, \cF_2\sub\cU_2,\dots$
such that the sets $\seq{\Un\cF}$ form an open cover of $X$.
A property introduced by Menger in~\cite{Menger24} was proved equivalent to this
covering property by Hurewicz in~\cite{Hure25}.
Thus, every compact space has Menger's property, and every space with Menger's property
is a Lindel\"of space, that is, one where every open cover has a countable subcover.

Every compact space is a Menger space, and every countable union of Menger spaces is Menger.
However, even among subsets of the real line there are large families of Menger spaces that are substantially
different from countable unions of compact spaces (e.g.,~\cite{MHP, sfh}).
Menger's property, which is central in the recent theory of \emph{selection principles}
(see~\cite{SaSchRP} and references therein), found
applications to seemingly unrelated notions in set theoretic and general topology and in
real analysis.

\brem
We mention two examples illustrating the importance of Menger's property in general and set theoretic
topology. This remark is independent of the remainder of the present paper, and we
refer the interested reader to the cited references
for definitions.

One of the major problems in set theoretic topology asks whether every
regular Lindel\"of space is a D-space. In the realm of Hausdorff spaces, the problem was
answered in the negative~\cite{SouSzep}. It turned out that
all \emph{Menger} spaces are D-spaces~\cite{AurD}.
Menger's property is still the most general natural class of spaces for
which a positive answer to the D-space problem is known.

In a series of papers (see~\cite{CanjarI, CanjarII} and references therein), a number of authors have
studied an important type of filters with a property introduced by Canjar. This property is related to the
theory of forcing: A filter has Canjar's property if the Mathias forcing notion associated to the filter does not add dominating reals.
It turned out that a filter has Canjar's property if and only if it is Menger in the standard, Cantor space topology~\cite{ChoZdo}.
This made a wide body of knowledge on Menger's property applicable to Canjar filters. In particular,
a number of earlier results follow immediately from this characterization.
\erem

Following Hurewicz~\cite{Hure25}, we restrict Menger's property to \emph{countable} open covers.
For Lindel\"of spaces, the two variations of Menger's property coincide, but otherwise the results obtained
are more general. This will be of importance to some applications at the end of this paper.

\bdfn
Let $X$ be a topological space.
A countable family $\cU$ of subsets of $X$ is an \emph{ascending cover} of $X$
if it is a cover of $X$ and there is an enumeration $\cU=\{\seq{V}\}$ such that
$V_1\subsetneq V_2\subsetneq\cdots$. Let $\Asc(X)$ be the family of open
covers of $X$ that \emph{contain} an ascending cover of $X$.
\edfn

We consider the family $P(X)$ of subsets of a set $X$
as a semigroup with the addition operator $\cup$.
Thus, for a family of sets $\cU\sub P(X)$, the set
$\FS(\cU)$ is comprised of all finite unions of members of $\cU$.
Only covers with no finite subcover constitute a challenge to
Menger's property.

\blem\label{lem:asc}
Let $X$ be a topological space. For each countable open cover $\cU$ with no finite subcover,
we have that $\FS(\cU)\in\Asc(X)$.\qed
\elem

For a topological space $X$, let $\Op(X)$ be the family of countable open covers of $X$.
A cover of $X$ is \emph{point-infinite} if every point of the space $X$ is contained in
infinitely many members of the cover.
Let $\Lambda(X)$ be the family of countable open point-infinite covers of $X$.
The proof of~\cite[Corollary~6]{coc1} establishes, in fact, that
$\sone(\Asc(X),\allowbreak \Lambda(X))$ holds whenever
$\sone(\Asc(X),\allowbreak \Op(X))$ does.

\bcor[Folklore]\label{cor:folk}
A topological space $X$ is Menger if and only if $\sone(\Asc(X),\allowbreak
\Lambda(X))$ holds.\qed
\ecor

Using a game theoretic theorem of Hurewicz, Scheepers proved in~\cite{OpPar} that
a space $X$ is Menger if and only if Alice does not have a winning strategy in the game
$\gfin(\Lambda(X),\Lambda(X))$, a variation of $\gone(\Lambda(X),\Lambda(X))$ where
Bob is allowed to choose any finite number of elements in each turn. Scheepers's theorem
is used in the following proof.

\bprp\label{prp:game}
A topological space $X$ is Menger if and only if Alice does not have a winning
strategy in the game $\gone(\Asc(X),\Lambda(X))$.
\eprp
\bpf
$(\Leftarrow)$ If Alice does not have a winning strategy in the game $\gone(\Asc(X),\Lambda(X))$,
then $\sone(\Asc(X),\Lambda(X))$ holds. Then $X$ is a Menger space.

$(\Rightarrow)$ Assume that Alice has a winning strategy in the game $\gone(\Asc(X),\Lambda(X))$.
Using this strategy, define a strategy for Alice in the game $\gfin(\Asc(X),\Lambda(X))$, as follows.
In the $n$-th inning, Alice's strategy proposes a cover containing an ascending cover.
Alice thins out this cover to make it ascending, and then removes from it the finitely many
elements chosen by Bob in the earlier innings.
This can only make Bob's task harder.
If Bob picks a finite subset $\cF_n$ of this ascending cover, Alice takes the largest set chosen by Bob,
$B_n$, and applies her original strategy, pretending that Bob chose only this set.

Assume that Bob won a play $(\cU_1,\cF_1,\cU_2,\cF_2,\dots)$ of the game
$\gfin(\Asc(X),\Lambda(X))$. Then $\Un_n\cF_n$ is a point-infinite cover of $X$.
Since the sets $\cF_n$ are disjoint, the set $\{\seq{B}\}$ is also a point-infinite cover of $X$,
and we obtain a play in the game $\gone(\Asc(X),\Lambda(X))$ that is won by Bob; a contradiction.
Thus, Alice has a winning strategy in the game $\gfin(\Asc(X),\Lambda(X))$.
Since $\Asc(X)\sub\Lambda(X)$, Alice has a winning strategy in the game
$\gfin(\Lambda(X),\Lambda(X))$. By Scheepers's Theorem, the space $X$ is not Menger.
\epf

With results proved thus far, we are ready to prove our main theorem.

\bthm\label{thm:MengerMT}
Let $(X,\tau)$ be a Menger space, and $\cU_1\sps\cU_2\sps\cdots$ be countable point-infinite
open covers of $X$ with no finite subcover.
For each finite coloring of the sets $\tau$ and $[\tau]^2$,
there are nonempty disjoint finite sets $\cF_1\sub\cU_1,\cF_2\sub\cU_2,\dots$
such that the sets $V_n := \Un\cF_n$, for $n\in\bbN$, have the following properties:
\be
\item The family $\{\seq{V}\}$ is a point-infinite cover of $X$.
\item The sets $\Un_{n\in F}V_n$ and $\Un_{n\in H}V_n$, for nonempty finite sets
$F<H$, are distinct.
\item All sets $\Un_{n\in F}V_n$, for nonempty finite sets $F\sub\bbN$, have the same color.
\item All sets $\{\,\Un_{n\in F}V_n, \Un_{n\in H}V_n\,\}$,
for nonempty finite sets $F<H$, have the same color.
\ee
Moreover, if $\cU_1=\{\seq{U}\}$, we may request that
the sets $F_n:=\set{m}{U_m\in\cF_n}$ satisfy
$F_1<F_2<\cdots$.
\ethm
\bpf
Enumerate $\cU_1=\{\seq{U}\}$. Consider the semigroup $(\tau,\cup)$.
We will work inside its subsemigroup $S:=\FS(\seq{U})$.
Let
$$\cA:=\set{A\sub S}{A\in\Asc(X)}.$$
The family $\cA$ is a superfilter:
Since $\cU_1$ has no finite subcover, the sequence $U_1,U_1\cup U_2,\dots$ has an ascending
subsequence. Thus, $\{U_1,U_1\cup U_2,\dots\}\in\cA$.
If $A\cup B\in\cA$, then the set $A\cup B$ contains an ascending cover $V_1\subsetneq V_2\subsetneq\cdots$,
and $A$ or $B$ must contain a subsequence of $\seq{V}$. Thus, $A\in\cA$ or $B\in\cA$.
The superfilter $\cA$ is translation invariant.
In particular, the superfilter $\cA$ is idempotent.

For each $n$, using that the cover $\cU_1$ has no finite subcover, fix an element
$x_n\in X\sm \Un_{i=1}^nU_i$.
For each $n$, let
$$\cV_n:=\set{V\in\FS(\set{U_m\in\cU_n}{m\ge n})}{x_1,\dots,x_{n-1}\in V}\sub\FS(U_n,U_{n+1},\dots).$$
(Note that $\cV_1=S$.)
For each $n$, the set $\set{U_m\in\cU_n}{m\ge n}$, being a cofinite subset of
the point-infinite cover $\cU_n$, is a (point-infinite) cover of $X$.
Since $\cU_n$ has no finite subcover,
we have that $\cV_n\in\Asc(X)$.
In particular, the sets $\cV_n$ are infinite.
We have that $\cV_1\sps\cV_2\sps\cdots$, and $\bigcap_n\cV_n=\emptyset$.
For each $n$, $\cV_n$ is a subsemigroup of $S$.
Thus, the sequence $\seq{\cV}$ is a free idempotent chain in $\cA$.

By Proposition~\ref{prp:game}, Alice does not have a winning strategy
in the game $\gone(\cA,\Lambda(X))$.
By Theorem~\ref{thm:abs},
for each finite coloring of the sets $S$ and $[S]^2$,
there are elements $V_1\in \cV_1, V_2\in \cV_2,\dots$ such that:
\be
\item The set $\{\seq{V}\}$ is in $\Lambda(X)$.
\item The sequence $\seq{V}$ is a proper sumsequence of $\seq{U}$.
\item The set $\FS(\seq{V})$ is monochromatic.
\item The sum graph of $\seq{V}$ is monochromatic.
\ee
The last assertion in the theorem is clear from the proof of Theorem~\ref{thm:abs}.
\epf

The assumption in Theorem~\ref{thm:MengerMT} that the space is Menger is necessary.
It is proved in~\cite{coc1} that being a Menger space is equivalent to the following property:
For each descending sequence $\cU_1\sps\cU_2\sps\cdots$ of countable point-infinite
open covers of $X$ with no finite subcover,
there are nonempty finite sets $\cF_1\sub\cU_1,\cF_2\sub\cU_2,\dots$
such that the family $\set{\Un\cF_n}{n\in\bbN}$ is a cover of $X$.

The following example shows that the Milliken--Taylor Theorem, and thus Hindman's Theorem,
is an instance of Theorem~\ref{thm:MengerMT} where Menger's property is trivial: a countable,
discrete space.

\bexm\label{exm:impH}
Consider Theorem~\ref{thm:MTFin}.
Let $X$ be the set of all cofinite subsets of $\bbN$, with the discrete topology.
Since the space $X$ is countable, it is a Menger space.

For each $n$, let $O_n:=\set{A\in X}{n\in A}$. The family $\{\seq{O}\}$
is a point-infinite open cover of $X$ with no finite subcover. According to our
conventions, for a set $F\in\Fin$ we have that $O_F=\Un_{n\in F}O_n$.

Let $S:=\FS(\seq{O})$. Then $S$ is a semigroup, and
the map $\Fin\to S$ defined by $F\mapsto O_F$ is a semigroup isomorphism.
Thus, a finite coloring of the set $[\Fin]^2$ may be viewed as a finite coloring of the set $[S]^2$.
Let $F_1<F_2<\cdots$ be nonempty finite sets such that the sets $V_n:=O_{F_n}$
satisfy assertion (4) of Theorem~\ref{thm:MengerMT}, and the sets $F_n$ are as requested
in Theorem~\ref{thm:MTFin}.
\eexm

The deduction of the classic theorems in Example~\ref{exm:impH} uses
a twist: It would have been more natural to consider the cover of $\bbN$ by
singletons, but there are 2-colorings of $\Fin$ with no monochromatic cover of $\bbN$ by disjoint
finite sets.

According to Example~\ref{exm:impH},
the Milliken--Taylor (or Hindman) Theorem may be viewed as a theorem
about countable open covers of countable sets, and Theorem~\ref{thm:MengerMT}
may be viewed as an extension of these theorems from countable spaces to
Menger spaces of arbitrary cardinality.
It is illustrative to compare this interpretation with Fern\'andez Br\'eton's 
impossibility result \cite{FerBre}:
For every set $S$, there is a $2$-coloring of the semigroup $\op{Fin}(S)$ of finite subsets of
$S$ such that no uncountable subsemigroup of $\op{Fin}(S)$ is monochromatic.
This demonstrates that the improvement must be on the qualitative side. In our case,
we color a countable object induced by a countable cover,
and obtain a monochromatic cover of the space---a
nontrivial assertion when the covered space is uncountable.

\section{Richer covers}

Let $X$ be a topological space, and $\cA$ and $\cB$ be families of covers of $X$.
Let $\ufin(\cA,\cB)$ be the property that, for covers $\seq{\cU}\in\cA$ with no finite subcover,
there are finite sets $\cF_1\sub \cU_1,\cF_2\sub \cU_2,\dots$,
such that $\{\seq{\Un\cF}\}\in\cB$.

Menger's covering property is the same as $\ufin(\Op(X),\Op(X))$.
A number of important covering properties are of the form $\ufin(\Op(X),\cB)$.
Some examples are provided in the survey~\cite{SaSchRP} and in the references therein.
By Lemma~\ref{lem:asc}, we have the following observation.

\bprp\label{prp:ufineqs1}
Let $X$ be a topological space, and $\cB$ be a family of covers of $X$.
The assertions $\ufin(\Op(X),\cB)$ and $\sone(\Asc(X),\cB)$ are equivalent.\qed
\eprp

Let $\Om(X)$ be the family of open covers $\cU$ of $X$ such that $X\notin\cU$ and every finite subset of $X$
is contained in some member of the cover. This family, introduced by Gerlits and Nagy~\cite{GN},
is central to the study of local properties in functions spaces.
The property $\ufin(\Op(X),\Om(X))$ was first considered by Scheepers~\cite{coc1}.
By the request that $X$ does not belong to any member of $\Om(X)$, the members of $\Om(X)$ are infinite.
Moreover, the family $\Om(X)$ is a superfilter on the topology $\tau$ of $X$.
If a cover $\cU\in\Om(X)$ is finer than another open
cover $\cV$ (in the sense that every member of $\cU$ is contained in some member of $\cV$) with $X\notin\cV$,
then $\cV\in\Om(X)$.
To cover additional important cases, we generalize these properties.

\bdfn
Let $(X,\tau)$ be a topological space. A family $\cB$ of open covers of $X$ is \emph{regular} if
it has the following properties:
\be
\item Whenever $\cU\cup\cV\in\cB$, we have that $\cU\in\cB$ or $\cV\in\cB$.
\item For each cover $\cU\in\cB$ and each finite-to-one function $f\colon\cU\to\tau\sm\{X\}$
with $U\sub f(U)$ for all $U\in\cU$, the image of $f$ is in $\cB$.
\ee
\edfn

Most of the important families of rich covers are regular.

\bexm
Let $X$ be a topological space.
The family $\Om(X)$ is regular.
The family $\Lambda(X)$ satisfies the second, but not the first, regularity condition.
Let $\Ga(X)$ be the family of infinite open covers of $X$ such that each point in $X$
is contained in all but finitely many members of the cover.
The family $\Ga(X)$ is regular.
The property $\ufin(\Op(X),\Ga(X))$ was introduced by Hurewicz~\cite{Hure25}.
Another well-studied regular family, denoted $\Tau^*(X)$, was introduced in~\cite{tautau}.
\eexm

In the following proof, we use the following observation.
It extends, by induction, to any finite number of ascending covers.

\blem
Let $\{\seq{U}\}$ and $\{\seq{V}\}$ be ascending covers of a set $X$, enumerated as such.
Then the set $\{U_1\cap V_1, U_2\cap V_2,\dots\}$ is an ascending cover of $X$.
\qed
\elem

\bthm\label{thm:game}
Let $(X,\tau)$ be a topological space, and $\cB$ be a regular family of open covers of $X$.
The following assertions are equivalent:
\be
\item $\ufin(\Op(X),\cB)$.
\item $\sone(\Asc(X),\cB)$.
\item Alice does not have a winning strategy in the game $\gone(\Asc(X),\cB)$.
\item Alice does not have a winning strategy in the game associated to $\ufin(\Op(X),\cB)$.
\ee
\ethm
\bpf
Proposition~\ref{prp:ufineqs1} asserts the equivalence of (1) and (2).
It is immediate that (4) implies (1).

$(3)\Impl (4)$: Assume that Alice has a winning strategy in the game associated to
$\ufin(\Op(X),\allowbreak\cB)$. By the definition of the selection principle $\ufin(\cA,\cB)$,
Alice's covers must not have finite subcovers.
By taking finite unions, turn every cover in Alice's strategy
into an ascending one. This only restricts the possible moves of Bob, and turns them
into moves in the game $\gone(\Asc(X),\cB)$.
Thus, we obtain a winning strategy for Alice in the latter game.

$(2)\Impl(3)$:
Assume that Alice has a winning strategy in the game $\gone(\Asc(X),\cB)$.
We encode this strategy as follows.
Let $\cU=\{\seq{U}\}$ be Alice's first move.
For each choice $U_{m_1}$ of Bob, let $\cU_{m_1}=\{\seq{U^{m_1}}\}$ be Alice's next move.
For each choice $U^{m_1}_{m_2}$ of Bob, let $\cU_{m_1,m_2}=\{\seq{U^{m_1,m_2}}\}$ be Alice's next move,
etc.
Thus, we have for each sequence $m_1,m_2,\dots,m_k\in\bbN$ a cover
$\cU_{m_1,m_2,\dots,m_k}=\{\seq{U^{m_1,m_2,\dots,m_k}}\}\in\Asc(X)$.

Thinning out the covers Alice plays only restricts Bob's moves.
Thus, we may assume that Alice plays ascending covers, and that
every cover played by Alice does not contain any of the finitely many elements played
by Bob in the earlier innings.
For a natural number $n$,
let $\{1,\dots,n\}^{\le n}:=\Un_{i=0}^n\{1,\dots,n\}^i$,
the set of all sequences of length at most $n$ taking values in $\{1,\dots,n\}$,
where the only sequence in $\{1,\dots,n\}^0$ is the empty sequence $\varepsilon$.
We define $U^\varepsilon_m:=U_m$ for all $m$.
For each $n$, set
$$\cV_n:=\Bigl\{\,
\seq{\bigcap_{\sigma\in\{1,\dots,n\}^{\le n}} U^\sigma}
\,\Bigr\}.$$
Then $\cV_n$ is an ascending cover of $X$.
By the property $\sone(\Asc(X),\cB)$,
there are elements $V_1\in\cV_1, V_2\in\cV_2,\dots$ such that $\{\seq{V}\}\in\cB$.

The cover $\{\seq{V}\}$ refines the cover $\cU$. Since $\cU$ has no finite subcover,
the set $\{\seq{V}\}$ is infinite.
We construct two parallel plays,
\beqn
(\cU,U_{m_1},\cU_{m_1},U^{m_1}_{m_3},\cU_{m_1,m_3},\dots);\\
(\cU,U_{m_2},\cU_{m_2},U^{m_2}_{m_4},\cU_{m_2,m_4},\dots),
\eeqn
according to Alice's strategy.
We use that Alice's covers are ascending.
\be
\item Pick a natural number $m_1>1$ such that
$$V_1\sub U_{m_1}\in\cU,$$
and $\{V_2,\dots,V_{m_1}\}\sm\{V_1\}\neq\emptyset$.

\item Each of the sets $V_2,\dots,V_{m_1}$ is contained in some member of the cover $\cU$.
Pick a natural number $m_2>m_1$ such that $U_{m_2}\neq U_{m_1}$,
$$V_2\cup\cdots\cup V_{m_1}\sub U_{m_2}\in\cU,$$
and $\{V_{m_1+1},\dots,V_{m_2}\}\sm\{V_1,\dots,V_{m_1}\}\ne\emptyset$.

\item For $n=3,4,\dots$:
\be
\item If $n$ is odd: Each of the sets $V_{m_{n-2}+1},\dots,V_{m_{n-1}}$ is contained in
some member of the cover $\cU_{m_1,m_3,\dots,m_{n-2}}$.
Pick a natural number $m_{n}>m_{n-1}$ such that the set $U:=U^{m_1,m_3,\dots,m_{n-2}}_{m_{n}}$
is distinct from all sets picked earlier,
$$V_{m_{n-2}+1}\cup\cdots\cup V_{m_{n-1}}\sub
U\in\cU_{m_1,m_3,\dots,m_{n-2}},$$
and $\{V_{m_{n-1}+1},\dots,V_{m_{n}}\}\sm\{V_1,\dots,V_{m_{n-1}}\}\neq\emptyset$.

\item If $n$ is even: Each of the sets $V_{m_{n-2}+1},\dots,V_{m_{n-1}}$ is contained in
some member of the cover $\cU_{m_2,m_4,\dots,m_{n-2}}$.
Pick a natural number $m_{n}>m_{n-1}$
such that the set $U:=U^{m_2,m_4,\dots,m_{n-2}}_{m_{n}}$
is distinct from all sets picked earlier,
$$V_{m_{n-2}+1}\cup\cdots\cup V_{m_{n-1}}\sub
U\in\cU_{m_2,m_4,\dots,m_{n-2}},$$
and $\{V_{m_{n-1}+1},\dots,V_{m_{n}}\}\sm\{V_1,\dots,V_{m_{n-1}}\}\neq\emptyset$.
\ee
\ee
Define a function
$$f\colon\{\seq{V}\}\to\{U_{m_1},U_{m_2},U^{m_1}_{m_3},U^{m_2}_{m_4},\dots\}$$
as follows.
\be
\item Map $V_1$ to $U_{m_1}$.
\item Map each element of the set $\{V_2,\dots,V_{m_1}\}\sm\{V_1\}$ to $U_{m_2}$.
\item For $n=3,4,\dots$ map each element of the set
$\{V_{m_{n-2}+1},\dots,V_{m_{n-1}}\}\sm\{V_1,\dots,V_{m_{n-2}}\}$
to $U^{m_1,m_3,\dots,m_{n-2}}_{m_{n}}$ if $n$ is odd, and to
$U^{m_2,m_4,\dots,m_{n-2}}_{m_{n}}$ if $n$ is even.
\ee
The function $f$ is as needed in property (2) of regular families of covers, and it is surjective.
Since the family $\{\seq{V}\}$ is in $\cB$ and $\cB$ is regular,
the set
$\{U_{m_1},U_{m_2},U^{m_1}_{m_3},U^{m_2}_{m_4},\dots\}$ is in $\cB$.
By Property (1) of regular families of covers,
one of the families
$\{U_{m_1},U^{m_1}_{m_3},\dots\}$
or
$\{U_{m_2},U^{m_2}_{m_4},\dots\}$
is in $\cB$.
It follows that Bob wins one of these two games against Alice's winning strategy; a contradiction.
\epf

\bthm\label{thm:ufinMT}
Let $(X,\tau)$ be a topological space, and $\cB$ be a regular family of open covers of $X$
(e.g., $\Omega(X)$, $\Tau^*(X)$, or $\Ga(X)$). Assume that $\ufin(\Op(X),\cB)$ holds.
Let $\cU_1\sps\cU_2\sps\cdots$ be countable point-infinite open covers of $X$ with no finite subcover.
For each finite coloring of the sets $\tau$ and $[\tau]^2$,
there are nonempty disjoint finite sets $\cF_1\sub\cU_1,\cF_2\sub\cU_2,\dots$
such that the sets $V_n := \Un\cF_n$, for $n\in\bbN$, have the following properties:
\be
\item The family $\{\seq{V}\}$ is in $\cB$.
\item The sets $\Un_{n\in F}V_n$ and $\Un_{n\in H}V_n$, for nonempty finite sets
$F<H$, are distinct.
\item All sets $\Un_{n\in F}V_n$, for nonempty finite sets $F\sub\bbN$, have the same color.
\item All sets $\{\,\Un_{n\in F}V_n, \Un_{n\in H}V_n\,\}$,
for nonempty finite sets $F<H$, have the same color.
\ee
Moreover, if $\cU_1=\{\seq{U}\}$,
we may request that the sets $F_n:=\set{m}{U_m\in\cF_n}$ satisfy
$F_1<F_2<\cdots$.
\ethm
\bpf
The proof is identical to that of Theorem~\ref{thm:MengerMT}, replacing $\Lambda(X)$ by
$\cB$ and using Theorem~\ref{thm:game} instead of Proposition~\ref{prp:game}.
\epf

In all of our theorems, the converse implications also hold.

\bprp
Let $X$ be a topological space, and $\cB$ be a regular family of open covers of $X$.
Assume that, for each
descending sequence $\cU_1\sps\cU_2\sps\cdots$ of countable point-infinite
open covers of $X$ with no finite subcover,
there are nonempty disjoint finite sets $\cF_1\sub\cU_1,\cF_2\sub\cU_2,\dots$
such that $\{\seq{\Un\cF}\}\in\cB$.
Then $\ufin(\Op(X),\cB)$ holds.
\eprp
\bpf
Let $\gimel(\cB)$ be the family of open covers $\cU$ of $X$ with no finite subcover,
such that there are disjoint finite sets $\seq{\cF}\sub\cU$ with $\{\seq{\Un\cF}\}\in\cB$.
By the first regularity property of $\cB$, we have that $\Lambda(X)\sps\gimel(\cB)$,
and by the premise of the proposition, $\Lambda(X)\sub\gimel(\cB)$.
By Scheepers's theorem, quoted after the proof of Theorem~\ref{thm:MengerMT},
the space $X$ is Menger.
By Corollary~10 and Lemma~11 of~\cite{GlCovs}, $\ufin(\Op(X),\cB)$ holds.
\epf

\section{Covers by more general sets}

\subsection{Borel covers}
Consider the variation of Menger's property, where covers by \emph{Borel sets} are considered.
Here, the restriction to countable covers is necessary to make the property nontrivial.\footnote{Otherwise,
the space could be covered by singletons, and then being Menger would be the same as being
countable.}
This property has its own history and applications (see, e.g.,~\cite{cbc} and the papers
citing it). As a rule, the results known for Menger's property extend to its Borel version~\cite{cbc},
and thus Theorem~\ref{thm:MengerMT} and its consequences also hold with ``open'' replaced
by ``Borel''.
The same assertion holds for the Borel versions of the other covering properties considered above.

In addition to open or Borel, one may consider other types of sets. As long as these types are preserved by the
basic operations used in the proof (mainly, finite intersections), the results obtained here apply
to countable covers by sets of the considered type.

\subsection{A combinatorial theorem}

Order the set $\NN$ by coordinate-wise comparison: $f\le g$ if $f(n)\le g(n)$ for all $n$.
Let $\fd$ be the minimal cardinality of a \emph{dominating family} $D\sub\NN$, that is,
such that for each function $f\in\NN$ there is a function $g\in D$ such that $f\le g$.
It is known that $\aleph_1\le\fd\le 2^{\aleph_0}$, but it is consistent that the cardinal
$\fd$ is strictly greater than $\aleph_1$ (more details are available in~\cite{BlassHBK}).
Let $D\sub\NN$ be a dominating family of cardinality $\fd$.
Then the property $\ufin(\Op(D),\Op(D))$ fails \cite{coc2}.
On the other hand, since we consider countable covers only,
the property $\ufin(\Op(X),\Op(X))$ holds for spaces $X$ of cardinality smaller than $\fd$ \cite{coc2}.
Thus, thinking of a cardinal number $\kappa$ as a discrete space of cardinality $\kappa$,
the following assertions are equivalent:
\be
\item $\kappa<\fd$.
\item $\ufin(\Op(\kappa),\Om(\kappa))$ holds.
\item $\ufin(\Op(\kappa),\Op(\kappa))$ holds.
\ee
By Theorem~\ref{thm:ufinMT}, we have the following purely combinatorial result.

\bthm\label{thm:comb}
Let $\kappa$ be a cardinal number smaller than $\fd$.
Let $\cU_1\sps\cU_2\sps\cdots$ be countable point-infinite covers of $\kappa$ with no finite subcover.
For each finite coloring of the sets $P(\kappa)$ and $[P(\kappa)]^2$,
there are nonempty disjoint finite sets $\cF_1\sub\cU_1,\cF_2\sub\cU_2,\dots$
such that the sets $A_n := \Un\cF_n$, for $n\in\bbN$, have the following properties:
\be
\item Every finite subset of $\kappa$ is contained in some set $A_n$.
\item The sets $\Un_{n\in F}A_n$ and $\Un_{n\in H}A_n$, for nonempty finite sets
$F<H$, are distinct.
\item All sets $\Un_{n\in F}A_n$, for nonempty finite sets $F\sub\bbN$, have the same color.
\item All sets $\{\,\Un_{n\in F}A_n, \Un_{n\in H}A_n\,\}$,
for nonempty finite sets $F<H$, have the same color.
\ee
Moreover, if $\cU_1=\{\seq{B}\}$, we may request that
the sets $F_n:=\set{m}{B_m\in\cF_n}$ satisfy
$F_1<F_2<\cdots$.\qed
\ethm

\section{Comments}

\subsection{Higher dimensions}
Our theorems also hold in dimensions larger than $2$, with minor modifications
in the proofs. For a natural number $d$, let
$[S]^d$ be the family of all $d$-element subsets of $S$.
We state the $d$-dimensional versions of
Theorem~\ref{thm:abs} and Theorem~\ref{thm:MengerMT}.
For brevity, the last part of Theorem~\ref{thm:7.2} is omitted.

\bthm
Let $S$ be a semigroup, and $d$ be a natural number.
Let $\cA$ be an idempotent superfilter on $S$, and $\cB$ be a family of subsets of $S$
such that Alice does not have a winning strategy in the game $\gone(\cA,\cB)$.
Let $\seq{a}$ be a sequence in $S$, and $A_1\sps A_2\sps\cdots$ be a free
idempotent chain in $\cA$ with $A_n\sub\FS(a_n,a_{n+1},\dots)$ for all $n$.
For each finite coloring of the set $[S]^d$,
there are elements $b_1\in A_1, b_2\in A_2,\dots$ such that:
\be
\item The set $\{\seq{b}\}$ is in $\cB$.
\item The sequence $\seq{b}$ is a proper sumsequence of $\seq{a}$.
\item The set $\set{\{b_{F_1},\dots,b_{F_d}\}}{F_1,\dots,F_d\in\Fin, F_1<\cdots<F_d}$
is monochromatic.
\ee
\ethm

\bthm\label{thm:7.2}
Let $(X,\tau)$ be a Menger space, and $d$ be a natural number.
For each descending sequence $\cU_1\sps\cU_2\sps\cdots$ of countable point-infinite
open covers of $X$ with no finite subcover,
and each finite coloring of the set $[\tau]^d$,
there are nonempty disjoint finite sets $\cF_1\sub\cU_1,\cF_2\sub\cU_2,\dots$
such that the sets $V_n := \Un\cF_n$, for $n\in\bbN$, have the following properties:
\be
\item The family $\{\seq{V}\}$ is a point-infinite cover of $X$.
\item The sets $\Un_{n\in F}V_n$ and $\Un_{n\in H}V_n$, for nonempty finite sets
$F<H$, are distinct.
\item All sets $\{\,\Un_{n\in F_1}V_n, \dots, \Un_{n\in F_d}V_n\,\}$,
for nonempty finite sets $F_1<\cdots<F_d$, have the same color.
\ee
\ethm

The $d$-dimensional versions of Theorems~\ref{thm:ufinMT} and~\ref{thm:comb} are similar.

\subsection{Proper sequences}
We have taken the approach of proper sequences, or having proper sumsequences,
to avoid pathological cases in theorems of Milliken--Taylor type.
Hindman and Strauss propose an unconditional approach in~\cite{HS}.
Corollary 18.9 in~\cite{HS} allows loops in the sum graph and considers colorings
of the set $[S]^1\cup[S]^2$.
In a manner similar to the proof of Proposition~\ref{prp:proper}, we obtain the
following observation.

\bprp\label{prp:sing}
Let $S$ be a semigroup, and consider the coloring $\chi$ of $[S]^1\cup [S]^2$ defined by
$\chi(\{a,b\}):=\card{\{a,b\}}$.
If a sequence $\seq{a}\in S$ has no proper sumsequence,
then every monochromatic sum graph of a sumsequence of $\seq{a}$ is a
singleton.\qed
\eprp

Since we may assume that any given finite coloring of the set $[S]^1\cup[S]^2$ is finer than
the one of Proposition~\ref{prp:sing}, there is no advantage
in this approach over that of Theorem~\ref{thm:MTS}.

\subsection{New covering properties}
Our results suggest a number of new covering properties that were not considered thus far,
and it remains unclear how exactly these relate to the classic ones. For example, the property
in Theorem~\ref{thm:MengerMT}, in the case where $\cU_n=\cU$ for all $n$, is formally weaker
than Menger's property. Is it equivalent to it?

\subsection{Additional directions}
Using the selection principle $\sfin$ and its corresponding game,
one obtains an abstract version of a theorem of Deuber and Hindman~\cite{DH},
and stronger forms of this theorem, in the spirit of the main theorem in Bergelson and
Hindman~\cite{BH89}. This direction will be pursued in a later project.

\ed